\documentclass[12pt]{amsart}
\usepackage{latexsym}
\usepackage{amsthm}
\usepackage{amsmath}
\usepackage{amsfonts}
\usepackage{amssymb}
\usepackage[dvips]{graphicx}
\input xy
\xyoption{all}

\newtheorem{theo}{Theorem}[section]
\newtheorem{lemma}[theo]{Lemma}
\newtheorem{propo}[theo]{Proposition}
\newtheorem{defi}[theo]{Definition}
\newtheorem{coro}[theo]{Corollary}
\newtheorem{rem}[theo]{Remark}

\newtheorem{exam}[theo]{Example}

\newcommand\hocolim{\operatorname{hocolim}}
\newcommand\holim{\operatorname{holim}}
\newcommand\BD{\boldsymbol{\Delta}}
\newcommand\Ind{\operatorname{Ind}}
\newcommand\HInd{\operatorname{HInd}}
\newcommand\Sind{\operatorname{Sind}}
\newcommand\HSind{\operatorname{HSind}}
\newcommand\HMod{\operatorname{HMod}}
\newcommand\Pre{\operatorname{Pre}}

\newcommand\Alg{\operatorname{Alg}}
\newcommand\HAlg{\operatorname{HAlg}}
\newcommand\op{\operatorname{op}}

\newcommand\SCat{\operatorname{\bf SCat}}

\newcommand\id{\operatorname{id}}
\newcommand\Id{\operatorname{Id}}
\newcommand\map{\operatorname{map}}

\newcommand\Set{\operatorname{\bf Set}}
\newcommand\SSet{\operatorname{\bf SSet}}
\newcommand\Sp{\operatorname{\bf S}}

\newcommand\Ho{\operatorname{Ho}}
\newcommand\Int{\operatorname{Int}}

\newcommand\ca{\mathcal {A}}

\newcommand\cc{\mathcal {C}}
\newcommand\cd{\mathcal {D}}

\newcommand\ch{\mathcal {H}}
\newcommand\ci{\mathcal {I}}
\newcommand\cj{\mathcal {J}}
\newcommand\ck{\mathcal {K}}
\newcommand\cl{\mathcal {L}}
\newcommand\cm{\mathcal {M}}

\newcommand\ct{\mathcal {T}}

\newcommand\cx{\mathcal {X}}

\date{May 26, 2006}
 
\begin{document}
\title{On homotopy varieties}
\author[J. Rosick\'{y}]
{J. Rosick\'{y}$^*$}
\thanks{ $^*$ Supported by the Ministry of Education of the Czech Republic under the project
MSM 0021622409.}
\address{\newline J. Rosick\'{y}\newline
Department of Algebra and Geometry\newline
Masaryk University, Faculty of Sciences\newline
Jan\'{a}\v{c}kovo n\'{a}m. 2a, Brno, Czech Republic\newline
rosicky@math.muni.cz
}
\begin{abstract}
Given an algebraic theory $\ct$, a homotopy $\ct$-algebra is a simplicial set where all equations from $\ct$
hold up to homotopy. All homotopy $\ct$-algebras form a homotopy variety. We will give a characterization of
homotopy varieties analogous to the characterization of varieties.  
\end{abstract}
\subjclass{ }

\maketitle

\section{Introduction}
Algebraic theories were introduced by F. W. Lawvere (see \cite{L1} and also \cite{L2}) in order to provide
a convenient approach to study algebras in general categories. An \textit{algebraic theory} is a small category
$\ct$ with finite products. Having a category $\ck$ with finite products, a $\ct$-\textit{algebra} in $\ck$ is
a finite product preserving functor $\ct\to\ck$. Algebras in the ca\-te\-go\-ry $\Set$ of sets are usual
many-sorted algebras. Algebras in the category $\SSet$ of simplicial sets are called simplicial algebras and they can 
be also viewed as simplicial objects in the ca\-te\-go\-ry of  algebras in $\Set$. 
In homotopy theory, one often needs to consider algebras up to homotopy -- a homotopy $\ct$-algebra is a functor 
$A:\ct\to\SSet$ such that the canonical morphism
$$
A(X_1\times\dots\times X_n)\to A(X_1)\times\dots\times A(X_n)
$$ 
is a weak equivalence for each finite product $X_1\times\dots\times X_n$ in $\ct$. These homotopy algebras have been 
considered in recent papers \cite{B1}, \cite{B2} and \cite{Be} but the subject is much older (see, e.g., \cite{BV}, 
\cite{Bc}, \cite{Ma} or \cite{St}).

A category is called a \textit{variety} if it is equivalent to the category $\Alg(\ct)$ of all $\ct$-algebras
in $\Set$ for some algebraic theory $\ct$. There is a characterization of varieties 
proved by F. W. Lawvere in the single-sorted case which can be immediately extended to varieties in general 
(cf. \cite{AR}, 3.25). Recent papers \cite {AR1} and \cite{ALR} reformulated this characterization by using the concept 
of a sifted colimit. Sifted colimits generalize filtered colimits -- while a category $\cd$ is filtered if colimits over 
$\cd$ commute with finite limits in $\Set$, a category $\cd$ is \textit{sifted} if colimits over $\cd$ commute
with finite products in $\Set$. 

Homotopy algebras will lead to homotopy varieties. But it is natural to consider \textit{simplicial algebraic
theories}, i.e., small simplicial categories $\ct$ with finite products. Homotopy algebras are then 
simplicial functors $\ct\to\SSet$ preserving finite products up to a weak equivalence. A ca\-te\-go\-ry will be called a 
\textit{homotopy variety} if it is homotopy equivalent to the category $\HAlg(\ct)$ of all homotopy $\ct$-algebras 
for some simplicial algebraic theory $\ct$.  Our main result is a characterization of homotopy varieties analogous to 
the just mentioned characterization of varieties. It uses the concept of a homotopy sifted homotopy colimit -- a category 
$\cd$ is \textit{homotopy sifted} if homotopy colimits over $\cd$ commute with finite products in $\SSet$. Homotopy 
sifted categories coincide with totally coaspherical categories in the sense of \cite{M}.  
  
Homotopy algebra is a part of a broader subject of a homotopy coherent category theory (see \cite{CP}).
More precisely, while usual algebra form a part of the theory of accessible and locally presentable categories
(see \cite{AR}), homotopy algebra belongs to the recently created theory of homotopy accessible and homotopy locally 
presentable categories (see \cite{Lu}, \cite{S}, \cite{TV1} and \cite{TV2}). J. Lurie \cite{Lu}
introduced homotopy accessible and homotopy locally presentable categories under the name of accessible 
$\infty$-categories and presentable $\infty$-categories. He works with CW-complexes instead of simplicial sets
and homotopy coherent functors in the sense of \cite{CP} and obtains a homotopy Giraud theorem characterizing
homotopy Grothendieck toposes. C. Simpson \cite{S} introduced a generalization of homotopy 
locally presentable categories using the language of Segal categories and called them, a little bit unfortunately, 
$\infty$-pretoposes. He characterized them as categories of fibrant and cofibrant objects of cofibrantly generated 
model categories. By D. Dugger \cite{D2}, homotopy locally presentable categories correspond in this way to combinatorial
model categories (i.e., cofibrantly generated and locally presentable). B. To\"{e}n and G. Vezzosi used 
the language of Segal categories to deal with homotopy Grothendieck toposes.

While presheaf categories $\Set^{\cc^{\op}}$, where $\cc$ is a small category, play a crucial role in the theory 
of accessible categories, homotopy accessible categories depend on the category $\SSet^{\cc^{\op}}$ of simplicial presheaves 
over a small simplicial category $\cc$. Recall that a simplicial category is a category enriched over $\SSet$. The category 
$\SSet^{\cc^{\op}}$ is equipped with the Bousfield-Kan (= projective) model category structure where both weak equivalences
and fibrations are pointwise. Then the homotopy presheaf category $\Pre(\cc)$ has to be taken as the full subcategory of 
$\SSet^{\cc^{\op}}$ consisting of objects which are both cofibrant and fibrant. It is the free completion of $\cc$ under 
homotopy colimits, which is analogous to $\Set^{\cc^{\op}}$ being the free completion of $\cc$ under colimits.  
In particular, the analogy of $\Set$ is the simplicial category $\Sp$ of fibrant simplicial sets, i.e., of Kan complexes.
All our simplicial categories have to be \textit{fibrant} in the sense that their hom-objects are Kan complexes. This makes
possible to define homotopy colimits in them -- just by using homotopy colimits in $\Sp$. Our framework is closely related to 
that based on $\infty$-categories or Segal categories. We expect that our results can be later placed to the context of 
quasi-categories which is under creation by A. Joyal (see \cite{J1} and \cite{J2}). 
 
\section{Simplicial categories}
A \textit{simplicial category} $\ck$ is a category enriched over the category $\SSet$ of simplicial sets. This means
that hom-objects $\hom(K,L)$ are simplicial sets and the composition of morphisms is given by simplicial maps. Morphisms
between simplicial categories are \textit{simplicial functors}, i.e., $F:\ck\to\cl$ is given by simplicial maps
$$
F_{K,L}:\hom(K,L)\to \hom(FK,FL)
$$ 
compatible with composition and unit. Morphisms between simplicial functors are \textit{simplicial natural 
transformations} $\varphi:F\to G$. They consist of morphisms $\varphi_K:FK\to GK$ for each $K$ in $\ck$ such that
the following diagram commutes for each pair of objects $K_1,K_2$ of $\ck$
$$
\xymatrix @=4pc{
\hom(K_1,K_2) \ar[r]^{F_{K_1,K_2}} \ar [d]_{G_{K_1,K_2}}&
\hom(FK_1,FK_2) \ar [d]^{\hom(FK_1,\varphi_{K_2})}\\
\hom(GK_1,GK_2) \ar [r]_{\hom(\varphi_{K_1},GK_2)} & \hom(FK_1,GK_2)
}
$$
(see \cite{H}, or \cite{Bo} for basic facts about enriched categories in general). 

We recall that an appropriate concept of (co)limits are weighted (co)limits. Given simplicial functors $D:\cd\to\ck$ and 
$G:\cd\to\SSet$, the limit $K$ of $D$ weighted by $G$ is defined by a simplicial natural isomorphism (in $X$)
$$
\hom(G,\hom(X,D))\cong \hom(X,K);
$$
hom's are always taken in appropriate simplicial categories. On the left side, it is the simplicial category of 
simplicial functors from $\cd$ to $\SSet$ (see \cite{Bo}); of course, $\hom(X,D):\cd\to\SSet$ is the composition
of $D$ and $\hom(X,-)$. Analogously, a colimit $K$ of $D:\cd\to\ck$ weighted by $G:\cd^{op}\to\SSet$ is given by a simplicial 
natural isomorphism 
$$
\hom(G,\hom(D,X))\cong \hom(K,X).
$$
 
Recall that a \textit{tensor} of a simplicial set $V$ and an object $K$ of a simplicial category $\ck$ is an object
$V\otimes K$ given by a simplicial natural isomorphism 
$$
\hom(V\otimes K,L)\cong\hom(V,\hom(K,L)).
$$
Dually, a \textit{cotensor} $K^V$ is given by
$$
\hom(L,K^V)\cong\hom(V,hom(L,K)).
$$
 
Model categories are taken in the sense of \cite{Ho} or \cite{H}. A \textit{simplicial 
model category} is a model category which is a simplicial category whose simplicial hom-sets are homotopically well behaved 
(see \cite{H} or \cite{GJ} for the precise definition). By \cite{D2}, every combinatorial model category is Quillen 
equivalent to a simplicial model category. Recall that a model category $\cm$ is called \textit{combinatorial} if the category
$\cm$ is locally presentable (cf. \cite{AR}) and its model structure is cofibrantly generated.  

There are well developed concepts of simplicial locally presentable ca\-te\-go\-ries and simplicial accessible categories 
(cf. \cite{K2}, \cite{BQ} and \cite{BQR}). Sim\-pli\-ci\-al locally presentable categories then correspond to weighted limit 
theories while simplicial accessible categories to theories specified by both weighted limits and weighted colimits. The desired
concepts of homotopy locally presentable categories and homotopy accessible categories should be based on homotopy limits
and homotopy colimits. In simplicial model categories, the definition of homotopy limits and homotopy colimits adopted
in \cite{H}, 18.1.8 and 18.1.1 make them a special case of weighted limits and weighted colimits (see \cite{H},
18.3.1); this observation goes back to \cite{BC}.
The corresponding weights form a homotopy invariant approximations of constant diagrams at a point.
The same definitions work in any simplicial category; in what follows, $B(\cx)$ denotes the nerve of the category $\cx$.

\begin{defi}\label{def2.1}
{\em
Let $\ck$ be a simplicial category, $\cd$ a small category and $D:\cd\to\ck$ a functor. Then the \textit{simplicial homotopy
colimit} $\hocolim_s D$ of $D$ is defined as the colimit of $D$ weighted by 
$$
B((-\downarrow\cd)^{op}):\cd^{op}\to\SSet.
$$

The \textit{simplicial homotopy limit} $\holim_s D$ of $D$ is defined as the limit of $D$ weighted by 
$$
B(\cd\downarrow-):\cd\to\SSet.
$$
}
\end{defi} 
 
Every simplicial category $\ck$ has a \textit{homotopy category} $\Ho(\ck)$; its objects are the same as that of $\ck$ and
$$
\hom_{\Ho(\ck)}(K,L)=\pi_0(\hom_\ck(K,L)),
$$
i.e., the set of morphisms from $K$ to $L$ in $\Ho(\ck)$ is the set of connected components of the simplicial set
of morphisms from $K$ to $L$ in $\ck$. Isomorphisms in $\Ho(\ck)$ are called \textit{homotopy equivalences} because
they coincide with usual homotopy equivalences for $\ck=\SSet$. We will use the notation $K\simeq L$ for homotopy
equivalent objects while $K\cong L$ will be kept for isomorphic objects.

The homotopy category of $\SSet$ in this sense is not the usual homotopy category of simplicial sets where isomorphisms are 
weak equivalences. In order to get the right one, one has to replace $\SSet$ by the simplicial category $\Sp$ of fibrant 
simplicial sets (i.e., of Kan complexes). Since homotopy equivalences coincide with weak equivalences here, simplicial 
$\Ho(\Sp)$ is the usual $\Ho(\SSet)$. 

Every simplicial functor $F:\ck\to\cl$ induces the functor 
$$\Ho(F):\Ho(\ck)\to \Ho(\cl).
$$ 
$F$ is called an \textit{equivalence} (see \cite{TV1}, 2.1.3) if

(1) the induced morphisms $\hom(K_1,K_2)\to \hom(F(K_1),F(K_2))$ are weak equivalences for all objects $K_1$
and $K_2$ of $\ck$ and

(2) each object $L$ of $\Ho(\cl)$ is isomorphic in $\Ho(\cl)$ to $\Ho(F)(K)$ for some object $K$ of $\ck$.

Let $\SCat$ denote the category of small simplicial categories and sim\-pli\-ci\-al functors. By \cite{Be1}, there is
a model category structure on $\SCat$ whose weak equivalences are the just defined equivalences. 
Fibrations are simplicial functors $F:\cc\to\cd$ satisfying two conditions (F1) and (F2) where the first one says 
that the simplicial maps $\hom(A,B)\to \hom(FA,FB)$ are fibrations of simplicial sets. In the special case when $\cd$ 
is the terminal simplicial category, (F1) says that hom-sets $\hom(A,B)$ are fibrant simplicial sets. Since (F2) is
automatic in this case, a small simplicial category $\cc$ is fibrant iff it has all $\hom(A,B)$ fibrant.

\begin{defi}\label{def2.2}
{\em
A simplicial category $\ck$ will be called \textit{fibrant} if all its hom-objects $\hom(A,B)$ are fibrant
simplicial sets.
}
\end{defi} 

For a simplicial model category $\cm$, $\Int(\cm)$ will denote its full subcategory consisting of objects
which are both cofibrant and fibrant. Then $\Int(\cm)$ is a fibrant simplicial category and its homotopy
category $\Ho(\Int(\cm))$ in the simplicial sense is equivalent to $\Ho(\cm)$ in the model category sense
(see \cite{H}). Recall that we have denoted $\Int(\SSet)$ by $\Sp$. Fibrant simplicial categories coincide
with categories enriched over $\Sp$.

The category $\Sp$ is closed in $\SSet$ under simplicial homotopy limits and under coproducts but it is not closed under
simplicial homotopy colimits in general. In order to get a concept of a homotopy colimit appropriate for $\Sp$, we have
to apply a fibrant replacement functor $R_f:\SSet\to\Sp$ to the simplicial homotopy colimit. We will call this new homotopy 
colimit \textit{fibrant} and denote it by $\hocolim_f$. Hence,
given a diagram $D:\cd\to\Sp$, we have
$$
\hocolim_f D = R_f(\hocolim_s D).
$$
The definition does not depend on a choice of a fibrant replacement functor because the resulting fibrant
homotopy colimits are always homotopy equivalent. From the model category point of view, there is no difference
between $\hocolim_s D$ and $\hocolim_f D$ because both objects are weakly equivalent.

Let $\cm$ be an arbitrary simplicial model category and consider a diagram $D:\cd\to\Int(\cm)$. We define its
\textit{fibrant homotopy colimit} $\hocolim_f D$ as $R_f(\hocolim_s D)$ where $R_f$ is a fibrant replacement
functor in $\cm$. Since $\hocolim_s D$ is cofibrant (see \cite{H}, 18.5.2), its fibrant replacement is both 
fibrant and cofibrant. Analogously, we define a \textit{fibrant homotopy limit} $\holim_f D$ as a cofibrant
replacement $R_c(\holim_s D)$. Since contravariant hom-functors of fibrant objects preserve weak equivalences
between cofibrant objects (see \cite{H}, 9.3.3), the simplicial sets $\hom(\hocolim_f D,A)$ and $\hom(\hocolim_s D,A)$
are weakly equivalent. Since covariant hom-functors of cofibrant objects preserve weak equivalences
between fibrant objects, the both simplicial sets are fibrant and thus they are homotopy equivalent. We get
$$
\hom(\hocolim_f D,A)\simeq\hom(\hocolim_s D,A)\cong\holim_s\hom(D,A).
$$
Analogously we get the formula
$$
\hom(A,\holim_f D)\simeq\holim_s\hom(A,D).
$$

\begin{defi}\label{def2.3}
{\em
Let $\ck$ be a fibrant simplicial category, $\cd$ a category and consider a diagram $D:\cd\to\ck$. We say that 
$\holim_f D$ is a \textit{fibrant homotopy limit} of $D$ if there are homotopy equivalences
$$
\delta_A:\hom(A,\holim_f D)\to\holim_s\hom(A,D)
$$
which are simplicially natural in $A$. 

Analogously, we define \textit{fibrant homotopy colimit} $\hocolim_f D$ of $D$ by the existence of homotopy
equivalences
$$
\delta_A:\hom(\hocolim_f D,A)\to\holim_s\hom(D,A).
$$
which are simplicially natural in $A$. 
}
\end{defi} 

In particular, we have the formulas
$$
\hom(A,\holim_f D)\simeq\holim_s\hom(A,D)
$$
and
$$
\hom(\hocolim_f D,A)\simeq\holim_s\hom(D,A).
$$
We will see in \ref{re3.1}(a) that $\holim_f D$ is determined uniquely up to a homotopy equivalence. In the case when 
$\ck=\Int(\cm)$ for a simplicial model category $\cm$, this definition coincides with the previous one. 

\begin{rem}\label{re2.4}
{\em
By the enriched Yoneda lemma, the simplicial natural transformation $\delta$ in the definition of the fibrant homotopy limit
is uniquely determined by its component
$$
\delta_{\holim_f D}:\hom(\holim_f D,\holim_f D)\to\holim_s\hom(\holim_f D,D)
$$
which uniquely corresponds to the morphism
$$
B(\cd\downarrow-)\to\hom(\hom(\holim_f D,\holim_f D),\hom(\holim_f D,D)).
$$
Since the codomain of this morphism is isomorphic to $\hom(\holim_f D,D)$, the simplicial natural transformation $\delta$
uniquely corresponds to the morphism
$$
\tilde{\delta}:B(\cd\downarrow-)\to\hom(\holim_f D,D).
$$

Analogously, the simplicial natural transformation $\delta$ in the definition of the fibrant homotopy colimit
is uniquely determined by its component
$$
\delta_{\hocolim_f D}:\hom(\hocolim_f D,\hocolim_f D)\to\holim_s\hom(D,\hocolim_f D)
$$
which uniquely corresponds to the morphism
$$
B(\cd^{\op}\downarrow-)\to\hom(\hom(\hocolim_f D,\hocolim_f D),\hom(D,\hocolim_f D)).
$$
Since the domain of this morphism is isomorphic to $B(-\downarrow\cd)^{\op}$ (see \cite{H}, 14.7.3) and the
codomain is isomorphic to $\hom(D,\hocolim_f D)$, the simplicial natural transformation $\delta$
uniquely corresponds to the morphism
$$
\tilde{\delta}:B(-\downarrow\cd)^{\op}\to\hom(D,\hocolim_f D).
$$

Morphisms $\tilde{\delta}$ correspond to limit cones for usual limits (see \ref{re3.1}(b)). We will sometimes
denote fibrant homotopy limits as pairs $(\holim_f D,\tilde{\delta})$.
Analogously, for fibrant homotopy colimits, $\tilde{\delta}$ correspond to colimit cocones.
}
\end{rem}

\begin{defi}\label{def2.5}
{\em
Let $F:\ck\to\cl$ be a simplicial functor between fibrant simplicial categories. We say that $F$ \textit{preserves} 
the fibrant homotopy limit of a diagram $D:\cd\to\ck$ if $(F\holim_f D,F\tilde{\delta})$ is a fibrant
homotopy limit of $FD$.

Analogously we define the preservation of fibrant homotopy colimits.
}
\end{defi} 

\begin{defi}\label{def2.6}
{\em
Let $G:\ck\to\cl$ and $F:\cl\to\ck$ be simplicial functors between fibrant simplicial categories. We say that $F$ is 
\textit{homotopy left adjoint} to $G$ if there are morphisms  
$$
\varphi_{K,L}:\hom(L,GK)\to\hom(FL,K) 
$$
and
$$
\psi_{K,L}:\hom(FL,K)\to\hom(L,GK)
$$
which are simplicially natural in $K$ and $L$ and such that $\psi_{K,L}$ is a homotopy inverse to $\varphi_{K,L}$ for 
each $K$ in $\ck$ and $L$ in $\cl$.
}
\end{defi} 

It implies that the induced functor $\Ho(F)$ is left adjoint to $\Ho(G)$. 
 
\section{Prestacks}

Let $\cc$ be a small simplicial category and consider the simplicial category $\SSet^{\cc^{\op}}$ of simplicial functors 
$\cc^{\op}\to\SSet$. We have the Yoneda embedding 
$$
Y_\cc:\cc\to\SSet^{\cc^{\op}}
$$ 
given by $Y(C)=\hom(-,C)$. The category $\SSet^{\cc^{\op}}$ has all weighted colimits and all weighted limits and
the Yoneda embedding $Y_\cc$ makes it the free completion of $\cc$ under weighted colimits. It also preserves all existing 
weighted limits (cf. \cite{K1}). Dually,
$$
\overline{Y}_\cc=Y^{\op}_{\cc^{\op}}:\cc\to (\SSet^\cc)^{\op}
$$ 
is the free completion of $\cc$ under weighted limits and preserves all existing weighted colimits. These free completions
exist for an arbitrary simplicial category -- one has to take \textit{small} simplicial functors into $\SSet$, i.e.,
small weighted (co)limits of hom-functors (see \cite{DL}).

For a small simplicial category $\cc$, $\SSet^{\cc^{\op}}$ is a simplicial combinatorial model category with respect to 
the projective (= Bousfield-Kan) model category structure. It means that weak equivalences and fibrations are pointwise.
(Trivial) cofibrations are then described as follows. We have the evaluation functors $E_C:\SSet^{\cc^{\op}}\to\SSet$,
$C\in\cc$; $E_C(F)=F(C)$. They are precisely the hom-functors 
$$
E_C=\hom(\hom(-,C),-).
$$ 
Each evaluation functor $E_C$ has a simplicial left adjoint 
$$
F_C=-\otimes\hom(-,C).
$$  
Now, cofibrations are cofibrantly generated by images in $F_C$, $C\in\cc$, of (generating) cofibrations in $\SSet$ and 
the same for trivial cofibrations. This procedure is described in \cite{H}, 11.6.1, for an ordinary category $\cc$ and 
\cite{CD} extends it to the simplicial category of small simplicial functors $\cc^{\op}\to\SSet$ for an arbitrary simplicial 
category $\cc$. The consequence is that all hom-functors $\hom(-,C)$ are cofibrant.

\begin{rem}\label{re3.1}
{\em
(a) Let $\ck$ be a fibrant simplicial category and assume that $L_1$ and $L_2$ are fibrant homotopy limits of a diagram 
$D:\cd\to\ck$. Let $\ck_0$ be a full subcategory of $\ck$ containing both $L_1$ and $L_2$. Then the hom-functors $\hom(-,L_1)$ 
and $\hom(-,L_2)$ are weakly equivalent in the projective model category $\SSet^{\ck_0^{\op}}$ with the functor 
$\holim_s\hom(-,D)$ restricted on $\ck_0$. Since they are cofibrant and fibrant, they are homotopy equivalent and thus $L_1$ 
and $L_2$ are homotopy equivalent.

(b) More generally, assume that we have an object $K$ in $\ck$ and a morphism 
$$
k:B(\cd\downarrow -)\to\hom(K,D).
$$
In the same way as in \ref{re2.4}, $k=\tilde\alpha$ for a simplicial natural transformation
$$
\alpha:\hom(-,K)\to\holim_s\hom(-,D).
$$
Let $\ck_0$ be a full subcategory of $\ck$ containing both $\holim_fD$ and $K$. Let 
$$
\gamma:R_cH\to H
$$
be a cofibrant replacement in $\SSet^{\ck_0^{\op}}$ of the restriction $H$ of the functor $\hom(-,\holim_fD)$ to $\ck_0$. 
Since hom-functors are cofibrant, there are simplicial natural transformations
$$
\delta':\hom(-,\holim_f D)\to R_cH
$$
and
$$
\alpha':\hom(-,K)\to R_cH 
$$
such that $\delta=\gamma\cdot\delta'$ and $\alpha=\gamma\cdot\gamma'$. Since $\gamma$ and $\delta$ are weak equivalences,
$\delta'$ is a weak equivalence and thus a homotopy equivalence because both $\hom(-,\holim_fD)$ and $R_cH$ are
cofibrant and fibrant. A homotopy inverse of $\delta'$ composed with $\alpha'$ gives a simplicial natural transformation
$$
\hom(-,K)\to\hom(-,\holim_fD)
$$ 
and thus a morphism $K\to\holim_fD$. This justifies our claim (cf. \ref{re2.4}) that $\tilde\delta$ plays the role of a limit 
cone. 
}
\end{rem}

\begin{theo}\label{th3.2}
Let $\cc$ be a small fibrant simplicial category. Then every object of $\SSet^{\cc^{\op}}$ is weakly equivalent to a simplicial
homotopy colimit of hom-functors.
\end{theo}
\begin{proof}
In the special case when $\cc$ is an ordinary category, the result was proved in \cite{D1}, 2.6. Let $\cc$ be an
arbitrary small simplicial category and consider its underlying ordinary category $\cc_0$. At first, we will assume
that $\cc$ has tensors $\Delta_n\otimes C$ for each $n$ and each $C$ in $\cc$.

Let $\BD$ be the category of non-empty finite ordinals and consider the product category $\cc_0\times\BD$ and the functor
$$
F:\cc_0\times\BD\to\cc
$$
given by the formula
$$
F(C,n)=\Delta_n\otimes C.
$$
Let
$$
F_\ast:\SSet^{\cc^{\op}}\to\SSet^{(\cc_0\times\BD)^{\op}}
$$
be the simplicial functor given by compositions with $F^{\op}$, i.e., $F_\ast(A)=A\cdot F^{\op}$. Then $F_\ast$ is the full 
embedding described in \cite{GJ}, p. 433. Moreover, it has a simplicial left adjoint
$$
F^\ast:\SSet^{(\cc_0\times\BD)^{\op}}\to\SSet^{\cc^{\op}}
$$
which is the weighted colimit preserving functor induced by the composition
$$
Y_\cc\cdot F:\cc_0\times\BD\to\SSet^{\cc^{\op}}.
$$
Thus $F^\ast$ preserves simplicial homotopy colimits and sends hom-functors to hom-functors. Since 
\begin{align*}
F^\ast(F_\ast(A))(C)&\cong\hom(\hom(-,C),(F^\ast\cdot F_\ast)(A))\\
&\cong\hom(\hom(-,(C,0)),F_\ast(A))\\
&\cong F_\ast(A)(C,0)\cong A(C), 
\end{align*}
$(F^\ast\cdot F_\ast)(A)$ is isomorphic to $A$ for each $A$ from $\SSet^{\cc^{\op}}$. Since the theorem is valid in 
$\SSet^{(\cc\times\BD)^{\op}}$, it suffices to show that $F^\ast$ preserves weak equivalences.

There is another simplicial functor
$$
\tilde{F}_\ast:\SSet^{\cc^{\op}}\to\SSet^{(\cc_0\times\BD)^{\op}}
$$
with a simplicial left adjoint
$$
\tilde{F}^\ast:\SSet^{(\cc_0\times\BD)^{\op}}\to\SSet^{\cc^{\op}}
$$
which preserves weak equivalences (see \cite{GJ}, Proposition IX.2.10). In the proof of this proposition, there is
found a pointwise homotopy equivalence
$$
\varrho:F_\ast\to\tilde{F}_\ast;
$$ 
it means that $\varrho_A:F_\ast(A)\to\tilde{F}_\ast(A)$ are homotopy equivalences for
each $A$ in $\SSet^{\cc^{\op}}$. The adjunction induces the morphism
$$
\sigma:\tilde{F}^\ast\to F^\ast
$$  
such that
$$
\hom(\sigma_B,A)\cdot\varphi_{A,B}=\tilde{\varphi}_{A,B}\cdot \hom(B,\varrho_A)
$$
where
$$
\varphi_{A,B}:\hom(B,F_\ast(A))\to \hom(F^\ast(B),A)
$$
and
$$
\tilde{\varphi}_{A,B}:\hom(B,\tilde{F}_\ast(A))\to \hom(\tilde{F}^\ast(B),A)
$$
denote the adjunction isomorphisms. Since hom-functors both preserve and reflect homotopy equivalences, $\sigma$
is a pointwise homotopy equivalence. Consequently, $F^\ast$ preserves weak equivalences.

Now, let $\cc$ be an arbitrary small fibrant simplicial category. Let
$$
G:\cc\to\hat\cc
$$
be the free completion of $\cc$ under tensors with $\Delta_n$, $n=1,2,\dots$. Then $\SSet^{\cc^{\op}}$ is
isomorphic to the full subcategory of $\SSet^{\hat\cc^{\op}}$ consisting of simplicial functors preserving
cotensors with $\Delta_n$, $n=1,2,\dots$. The inclusion of $\SSet^{\cc^{\op}}$ to $\SSet^{\hat\cc^{\op}}$ 
has a simplicial left adjoint 
$$
G_\ast:\SSet^{\hat\cc^{\op}}\to\SSet^{\cc^{\op}}
$$ 
given by compositions with $G$. Since $G_\ast$ preserves weak equivalences and the claim is valid in $\hat\cc$,
it suffices to show that images of hom-functors in $G_\ast$ are weakly equivalent to hom-functors. But we have
$$
G_\ast\hom(-,C^{\Delta_n})\cong \hom(-,C)^{\Delta_n}
$$
and the latter functor is weakly equivalent to $\hom(-,C)$ because $\Delta_n$ is contractible and $\cc$ is fibrant. 
\end{proof}

\begin{lemma}\label{lem3.3}
Let $G:\ck\to\cl$ be a simplicial functor between fibrant simplicial categories and $F:\cl\to\ck$ its homotopy left
adjoint. Then $F$ preserves fibrant homotopy colimits and $G$ preserves fibrant homotopy limits.
\end{lemma}
\begin{proof}
Let $D:\cd\to\ck$ be a diagram. We get simplicial natural transformations
$$
\hom(F\hocolim_f D,-)\to\hom(\hocolim_f D,G(-)),
$$
$$
hom(\hocolim_f D,G(-))\to\holim_s \hom(D,G(-))
$$
and
$$
\hom(\hocolim_f FD,-))\to\holim_s\hom(FD,-)
$$
whose components are homotopy equivalences. Since compatible weak equivalences between diagrams of fibrant objects
induce a weak equivalence of their simplicial homotopy limits, the functors 
$$
\hom(F\hocolim_fD,-)
$$ 
and
$$
\hom(\hocolim_fFD,-)
$$  
are weakly equivalent and thus homotopy equivalent. This implies that $F$ preserves fibrant 
homotopy colimits. The statement about $G$ is dual.
\end{proof}

\begin{defi}\label{def3.4}
{\em
Let $\cc$ be a small fibrant simplicial category. We put
$$
\Pre(\cc)=\Int(\SSet^{\cc^{\op}}).
$$
A simplicial functor $F:\cc^{\op}\to\SSet$ belonging to $\Pre(\cc)$ is called a \textit{prestack} on $\cc$. 
}
\end{defi} 

Prestacks of $\cc$ are precisely simplicial functors $\cc^{op}\to\Sp$ which are cofibrant objects in the Bousfield-Kan 
model category structure on $\SSet^{\cc^{\op}}$. Our terminology is in accordance with \cite{Lu} and \cite{TV2}.
Since $\cc$ is fibrant, all hom-functors $\hom(-,C)$ are prestacks because they are always cofibrant in 
$\SSet^{\cc^{\op}}$. Thus we get the Yoneda embedding
$$
Y_\cc:\cc\to\Pre(\cc).
$$
 
\begin{theo}\label{th3.5}
Let $\cc$ be a small fibrant simplicial category. Then $\Pre(\cc)$ is the free completion of $\cc$ under fibrant 
homotopy colimits.
\end{theo}
\begin{proof} The statement means that $\Pre(\cc)$ has fibrant homotopy colimits and for each simplicial functor
$F:\cc\to\ck$ with $\ck$ fibrant and having fibrant homotopy colimits there is, up to a homotopy equivalence, a unique 
simplicial functor 
$$
F^\ast:\Pre(\cc)\to\ck
$$ 
which preserves fibrant homotopy colimits and satisfies $F^\ast\cdot Y_\cc\simeq F$. Add that a homotopy equivalence 
of functors means a pointwise homotopy equivalence.

We know that $\Pre(\cc)$ has fibrant homotopy colimits and, by \ref{th3.2}, each prestack is homotopy equivalent
to a fibrant homotopy colimit of hom-functors. This yields a unique candidate, up to homotopy equivalence, of an
extension $F^\ast$ of $F$ which preserves fibrant homotopy colimits. We have to prove that $F^\ast$ is a simplicial 
functor which preserves all fibrant homotopy colimits. 
 
Since fibrant homotopy colimits in $\SSet^{\op}$ are fibrant homotopy limits in $\SSet$, each hom-functor
$$
\hom(-,K):\ck\to\SSet^{\op}
$$
preserves fibrant homotopy colimits. Thus our candidates satisfy
$$
\hom(-,K)\cdot F^\ast\simeq (\hom(-,K)\cdot F)^\ast
$$
for each object $K$ in $\ck$. Hence it suffices to prove our claim for $\ck=\Sp^{op}$. Let $F:\cc\to\Sp^{\op}$ be
a simplicial functor.

Since $\SSet^{\op}$ is a model category, there is a left Quillen functor
$$
G:\SSet^{\cc^{\op}}\to\SSet^{\op}
$$
with $G\cdot Y_\cc=F$ (see \cite{D1}). Since $G$ preserves simplicial homotopy colimits (up to a weak equivalence, see 
\cite{D1}) and weak equivalences of cofibrant objects, it preserves fibrant homotopy colimits. Therefore $F^\ast$ is homotopy 
equivalent to $R_f^{\op}\cdot G$, which proves the claim.
\end{proof}

\begin{rem}\label{re3.6}
{\em
Since fibrant homotopy colimits are pointwise, they are preserved by the evaluation functors $E_C:\Pre(\cc)\to\Sp$,
$C\in\cc$; $E_C(F)=F(C)$. Since the evaluation functors are precisely the hom-functors 
$$
E_C=\hom(\hom(-,C),-),
$$ 
representable functors $\hom(-,C)$ are \textit{homotopy absolutely presentable} in the sense that their hom-functors 
$\hom(\hom(-,C),-)$ preserve all fibrant homotopy colimits. We also have
$$
E_C\cdot Y_\cc = \hom(C,-).
$$  
}
\end{rem}

\begin{defi}\label{def3.7}
{\em
An object $K$ of a fibrant simplicial category $\ck$ is called \textit{homotopy finitely presentable} provided that its 
hom-functor 
$$
\hom(K,-):\ck\to\Sp
$$ 
preserves filtered fibrant homotopy colimits.
}
\end{defi} 

Recall that filtered homotopy colimits are over diagrams $D:\cd\to\ck$ where $\cd$ is a filtered category 
(cf. \cite{AR}). A homotopy colimit of a dia\-gram $\cd\to\ck$ is called \textit{finite} if the category $\cd$ has
finitely many morphisms.
 
\begin{propo}\label{prop3.8}
In $\Sp$, filtered fibrant homotopy colimits commute with finite fibrant homotopy limits.
\end{propo}
\begin{proof}
The statement means that, given a diagram $D:\ci\times\cj\to\Sp$ with $\ci$ filtered and $\cj$ finite, the canonical morphism 
$$
c:\underset{\ci}{\hocolim_f}\underset{\cj}{\holim_f} D(i,j) \to \underset{\cj}{\holim_f}\underset{\ci}{\hocolim_f} D(i,j)
$$
is a homotopy equivalence. By \cite{BK}, XII., 3.5(ii), filtered simplicial homotopy colimits are weakly equivalent 
to filtered colimits in $\SSet$. Since $\Sp$ is closed in $\SSet$ under filtered colimits, filtered fibrant homotopy colimits 
in $\Sp$ are homotopy equivalent with filtered colimits. Thus the result is a consequence of the fact that filtered colimits 
commute with finite weighted limits in $\SSet$ (see \cite{BQ}).  
\end{proof}

\begin{propo}\label{prop3.9}
Let $\ck$ be a fibrant simplicial category. Then a finite fibrant homotopy colimit of homotopy finitely
presentable objects is homotopy finitely presentable.
\end{propo}
\begin{proof}
Let $J:\cj\to\ck$ a finite diagram with homotopy finitely presentable values. We have to prove that 
$\hocolim_f J$ is homotopy finitely presentable. Let $I:\ci\to\ck$ be a filtered diagram. Then, by \ref{prop3.8}, 
we have  
\begin{align*}
\hom(\hocolim_f J,\hocolim_f I)&\simeq \holim_f\hom(J,\hocolim_f I)\\
&\simeq \underset{\cj}{\holim_f}\underset{\ci}{\hocolim_f} \hom(J,I)\\
&\simeq\underset{\ci}{\hocolim_f}\underset{\cj}{\holim_f} \hom(J,I)\\
&\simeq \hocolim_f \hom(\hocolim_f J,I).
\end{align*}
Thus $\hocolim_f J$ is homotopy finitely presentable. 
\end{proof}

\section{Homotopy varieties}
Consider a category $\ck$ with binary products and diagrams $D_1:\cd_1\to\ck$ and $D_2:\cd_2\to\ck$. We form the diagram
$$
D_1\times D_2:\cd_1\times\cd_2\to\ck
$$
by means of the formula $(D_1\times D_2)(d_1,d_2)=D_1d_1\times D_2d_2$ (do not confuse it with the product functor
$\cd_1\times\cd_2\to\ck\times\ck$). 

\begin{defi}\label{def4.1}
{\em
Let $\ck$ be a fibrant simplicial category having fibrant homotopy colimits and binary products. We say that
\textit{fibrant homotopy colimits distribute over binary products} in $\ck$ provided that
$$
\hocolim_f(D_1\times D_2)\simeq \hocolim_f D_1\times\hocolim_f D_2
$$
for every pair of diagrams $D_1:\cd_1\to\ck$ and $D_2:\cd_2\to\ck$. 
}
\end{defi}

\begin{propo}\label{prop4.2}
In $\Sp$, fibrant homotopy colimits distribute over binary products.
\end{propo}
\begin{proof}
Consider diagrams $D_1:\cd_1\to\Sp$ and $D_2:\cd_2\to\Sp$. Since the functor $-\times Y:\SSet\to\SSet$ has
the simplicial right adjoint $\hom(Y,-)$, it preserves simplicial homotopy colimits. Thus simplicial homotopy
colimits distribute over binary products in $\SSet$. Hence it suffices to know that a product $w_1\times w_2$
of weak equivalences $w_1$ and $w_2$ in $\SSet$ is a weak equivalence. This fact can be deduced as follows.

Since the geometric realization functor $|\,|:\SSet\to{\bf K}$ to the full subcategory $\bf K$ of the category 
of topological spaces consisting of compactly generated spaces preserves finite limits (see \cite{Ho}, 3.2.4), we have 
$|w_1\times w_2|=|w_1|\times|w_2|$. Now $w$ is a weak equivalence in $\SSet$ iff $|w|$ is 
a weak equivalence in $\bf K$ and the product of weak equivalences in $\bf K$ is a weak equivalence by \cite{H}, 18.5.3 
(because products are homotopy limits and all objects are fibrant in $\bf K$). 
\end{proof}

\begin{defi}\label{def4.3}
{\em
A small category $\cd$ will be called \textit{homotopy sifted} provided that fibrant homotopy colimits over $\cd$
commute in $\Sp$ with finite products.
}
\end{defi} 

Explicitly, $\cd$ is homotopy sifted iff it is nonempty (thus fibrant homotopy colimits over $\cd$ commute with 
the empty product) and, given diagrams $D_1,D_2:\cd\to\Sp$, then the canonical morphism
$$
\hocolim_f(D_1\otimes D_2)\to\hocolim_f D_1\times\hocolim_f D_2
$$
is a homotopy equivalence. Here, the diagram $D_1\otimes D_2:\cd\to\Sp$ is given by
$$
(D_1\otimes D_2)(d)=D_1d\times D_2d.
$$
In fact, $D_1\otimes D_2$ is the product of $D_1$ and $D_2$ in $\SSet^{\cd}$. 

The following theorem is analogous to the characterization of sifted colimits (see \cite{AR1} or \cite{ARV}). 
Recall that a functor $F:\ck\to\cl$ is called \textit{homotopy final} provided that for
every object $L$ of $\cl$ the comma-category $L\downarrow F$ is aspherical, i.e., its nerve $B(L\downarrow F)$
is weakly equivalent to the point (see \cite{H}, 19.6.1). Every homotopy final functor is final because
the latter means that all comma-categories $L\downarrow F$ are non-empty and connected.
 
\begin{theo}\label{th4.4}
A small category $\cd$ is homotopy sifted iff $\cd$ is nonempty and the diagonal functor 
$\Delta:\cd\to\cd\times\cd$ is homotopy final.
\end{theo}
\begin{proof}
Given diagrams $D_1,D_2:\cd\to\Sp$, we have
$$
D_1\otimes D_2=(D_1\times D_2)\Delta.
$$
By \cite{H}, 19.6.7 and 19.6.12, $\Delta$ is homotopy final iff the induced map 
$$
\hocolim_s DF\to\hocolim_s D
$$
is a weak equivalence for every diagram $D:\cd\times\cd\to\SSet$. This is clearly the same as
$$\hocolim_f DF\to\hocolim_f D
$$
being a homotopy equivalence for every diagram $D:\cd\times\cd\to\Sp$. Consequently, $\cd$ is homotopy sifted provided
that $\Delta$ is homotopy final. Conversely, since the proof of \cite{H}, 19.6.12 only uses functors
$$
D=\hom((d_1,d_2),-)=\hom(d_1,-)\times\hom(d_2,-),
$$
$\Delta$ is homotopy final whenever $\cd$ is homotopy sifted.
\end{proof}
 
\begin{rem}\label{re4.5}
{\em
(1) A category $\cd$ is homotopy sifted iff all comma-categories $(d_1,d_2)\downarrow\Delta$, where $d_1,d_2$
are objects from $\cd$, are aspherical. Hence $\cd$ is homotopy sifted iff $\cd^{\op}$ is totally aspherical 
in the sense of \cite{M}, 1.6.3.

(b) By \ref{prop3.8}, each filtered category is homotopy sifted. But it also follows from the fact that
every filtered category $\cd$ is aspherical because it is a filtered colimit of categories $d\downarrow\cd$
having the initial object (see \cite{Q}).

(c) Every category $\cd$ with finite coproducts is homotopy sifted (see \cite{M}, 7.4). It immediately follows
from the fact that $d_1\amalg d_2$ is the initial object in $(d_1,d_2)\downarrow\cd$.

(d) Every homotopy sifted category is sifted because $\Delta$ is final provided that it is homotopy final.

(e) Recall that a \textit{reflexive coequalizer} is defined a coequalizer of a pair of morphisms $h,k:A\to B$
which have a common section $m:B\to A$, i.e., such that $hm=km=\id_B$; such pairs are called \textit{reflexive}
(see \cite{ARV}). A \textit{ fibrant homotopy reflexive coequalizer} is defined as a fibrant homotopy coequalizer 
of a reflexive pair. Reflexive coequalizers form an important kind of sifted categories (see \cite{AR1} or 
\cite{ARV}). But they are not homotopy sifted -- a direct inspection shows that the comma category 
$(A,A)\downarrow\cd$ is not aspherical (it is connected but not 2-connected); $\cd$ denotes a reflexive pair.

(f) The reflexive pair is the full subcategory of the category $\BD^{\op}$ consisting of ordinals $1,2$.
The whole category $\BD^{\op}$ is homotopy sifted following \cite{M}, 1.6.13.
}
\end{rem}
 
\begin{defi}\label{def4.6}
{\em
An object $K$ of a fibrant simplicial category $\ck$ is called \textit{homotopy strongly finitely
presentable} provided that its hom-functor $hom(K,-):\ck\to\Sp$ preserves homotopy sifted fibrant homotopy colimits.
}
\end{defi} 

\begin{rem}\label{re4.7}
{\em
By \ref{re4.5}(b), every homotopy strongly finitely presentable object is homotopy finitely presentable.
}
\end{rem}

\begin{propo}\label{prop4.8}
A finite coproduct of homotopy strongly finitely presentable objects is homotopy strongly finitely
presentable.
\end{propo}
\begin{proof}
The proof is analogous to that of \ref{prop3.9}.
\end{proof}

\begin{propo}\label{prop4.9}
Let $G:\ck\to\cl$ be a simplicial functor between fibrant simplicial categories which has a homotopy left adjoint 
$F:\cl\to\ck$. Then $F$ preserves homotopy strongly finitely presentable objects provided that $G$ preserves 
homotopy sifted fibrant homotopy colimits.
\end{propo}
\begin{proof}
Assume that $G$ preserves homotopy sifted fibrant homotopy colimits. We have to show that for each homotopy strongly 
finitely presentable object $L$ of $\cl$ the object $FL$ is homotopy strongly finitely presentable as well.

Let $\cd$ be a homotopy sifted category and consider a diagram $D:\cd\to\ck$. We have
\begin{align*}
\hom(FL,\hocolim_f D)&\simeq \hom(L,G(\hocolim_f D))\\
&\simeq \hom(L,\hocolim_f GD)\\
&\simeq \hocolim_f hom(L,GD)\\
&\simeq \hocolim_f \hom(FL,D). 
\end{align*}
Hence $FL$ is homotopy strongly finitely presentable in $\ck$.
\end{proof}

\begin{defi}\label{def4.10}
{\em
A fibrant simplicial category $\ck$ will be called  a \textit{homotopy variety} provided that it has fibrant 
homotopy colimits and has a set $\ca$ of homotopy strongly finitely presentable objects such that every object 
of $\ck$ is a homotopy sifted fibrant homotopy colimit of objects from $\ca$. 
}
\end{defi} 

\begin{propo}\label{prop4.11}
Let $\cc$ be a small fibrant simplicial category. Then the category $\Pre(\cc)$ is a homotopy variety. 
\end{propo}
\begin{proof}
Let $\bar\cc$ be the closure of $Y(\cc)$ under finite coproducts in $\Pre(\cc)$. By \ref{re3.6} and \ref{prop4.8}, 
each object of $\bar\cc$ is homotopy strongly finitely presentable in $\Pre(\cc)$. For each object
$A$ in $\Pre(\cc)$, the comma-ca\-te\-go\-ry $\bar\cc\downarrow A$ has finite coproducts. By \ref{re4.5}(c),
$\bar\cc\downarrow A$ is homotopy sifted. Since $A$ is the fibrant homotopy colimit of the projection
$\bar\cc\downarrow A\to\Pre(\cc)$ (see \ref{th3.5}), the category $\Pre(\cc)$ is a homotopy variety. 
\end{proof}

\begin{defi}\label{def4.12}
{\em
A \textit{ simplicial algebraic theory} is defined as a small fibrant simplicial category $\ct$ having finite products.

A \textit{homotopy} $\ct$-\textit{algebra} is a simplicial functor $A:\ct\to\Sp$ belonging to $\Pre(\ct^{\op})$ such that 
the canonical morphism
$$
A(X_1\times\dots\times X_n)\to A(X_1)\times\dots\times A(X_n)
$$
is a homotopy equivalence for each finite product $X_1\times\dots\times X_n$ in $\ct$.

We will denote by $\HAlg(\ct)$ the full subcategory of $\Pre(\ct^{op})$ consisting of all homotopy
$\ct$-algebras.
}
\end{defi} 

\begin{exam}\label{ex4.13}
{\em
Let $\ct_0$ be the algebraic theory of one binary operation $m$. It means that $\ct_0$ has objects $X_0,X_1,\dots,
X_n,\dots$ and morphisms are generated by $m:X_2=X_1\times X_1\to X_1$. Then a $\ct_0$-algebra $A$ is a simplicial set
$A(X_1)$ equipped with a binary operation $A(m):A(X_1)\times A(X_1)\to A(X_1)$. Let $\ct_1$ be the simplicial 
algebraic theory obtained from $\ct_0$ by adding a one-dimensional simplex to $\hom(X_3,X_1)$ from the point
$m(m\times \id)$ to $m(\id\times m)$. It means that we have the corresponding simplicial map
$$
h:\Delta_1\to\hom(X_3,X_1).
$$
Given a $\ct_1$-algebra $A$, we get the composition 
$$
\Delta_1\to\hom(X_3,X_1)\to\hom(A(X_1)^3,A(X_1)
$$
of $h$ with $A_{X_3,X_1}$. This composition corresponds to the simplicial map
$$
\Delta_1\times A(X_3)\to A(X_1),
$$
which is a homotopy from $A(m)(A(m)\times\id)$ to $A(m)(\id\times A(m))$. In this way we can get strongly homotopy
associative algebras of \cite{St} as algebras for a suitable simplicial algebraic theory. Homomorphisms of
these algebras strictly preserve the multiplication.

On the other hand, if $\ct_2$ is the algebraic theory of one associative binary operation then homotopy
$\ct_2$-algebras are simplicial sets equipped with a homotopy associative multiplication and homomorphisms
preserve the operation up to homotopy.
}
\end{exam}

There is proved in \cite{B1} and \cite{Be} that, for a simplicial algebraic theory $\ct$, each homotopy
$\ct$-algebra is weakly equivalent to a strict $\ct$-algebra in a suitable model category structure on
$\SSet^{\ct}$.

\begin{propo}\label{prop4.14}
Let $\ct$ be a simplicial algebraic theory. Then the simplicial category $\HAlg(\ct)$ is closed in $\Pre(\ct^{\op})$ 
both under fibrant homotopy limits and homotopy sifted fibrant homotopy colimits.
\end{propo}
\begin{proof}
Consider a diagram $D:\cd\to\HAlg(\ct)$. We have
\begin{align*}
(\underset{d}{\holim_f} Dd)(X_1\times\dots\times X_n)&\simeq \underset{d}{\holim_f} Dd(X_1\times\dots\times X_n)\\
&\simeq \underset{d}{\holim_f} Dd(X_1)\times\dots\times \underset{d}{\holim_f} Dd(X_n).
\end{align*}
Thus $\HAlg(\ct)$ is closed in $\Pre(\ct^{\op})$ under fibrant homotopy limits. Since homotopy sifted fibrant homotopy
colimits commute in $\Sp$ with finite products, we analogously prove that $\HAlg(\ct)$ is closed in $\Pre(\ct^{\op})$ 
under homotopy sifted fibrant homotopy colimits.
\end{proof}

\begin{theo}\label{th4.15}
A fibrant simplicial category $\ck$ is a homotopy variety iff it is equivalent to $\HAlg(\ct)$ for some 
simplicial algebraic theory $\ct$. 
\end{theo}
\begin{proof}
I. Let $\ct$ be a simplicial algebraic theory. Consider a finite product diagram
$$
p_i:X_1\times\dots\times X_n\to X_i\quad i=1,\dots,n
$$
in $\ct$. Let
$$
m_{X_1\dots X_n}:\hom(X_1,-)\amalg\dots\amalg\hom(X_n,-)\to\hom(X_1\times\dots\times X_n,-)
$$
be the morphism induced by 
$$
\hom(p_i,-):\hom(X_i,-)\to\hom(X_1\times\dots\times X_n,-).
$$
Let $A:\cc\to\Sp$ be a functor belonging to $\Pre(\ct^{\op})$. Since 
$$
\hom(\hom(X_1\times\dots\times X_n,-),A)\cong A(X_1\times\dots\times X_n)
$$
and
$$
\hom(\hom(X_1,-)\amalg\dots\amalg\hom(X_n,-),A)\cong A(X_1)\times\dots\times A(X_n),
$$
the functor $A$ is a homotopy $\ct$-algebra iff $\hom(m_{X_1\dots X_n},A)$ is a homotopy equivalence for each finite
product diagram in $\ct$. 

Let $\cm$ be the set of all morphisms $m_{X_1\dots X_n}$. Recall that an object $A$ of $\SSet^\ct$ is \textit{homotopy
orthogonal} to $\cm$ if 
$$
\map(m_{X_1\dots X_n},A)
$$
is a weak equivalence for each $m_{X_1\dots X_n}$ from $\cm$ (see \cite{H}, 17.8.5). Here, $\map(B,A)$ denotes a homotopy 
function complex. Let $\cm^\bot$ be the full subcategory of $\SSet^\ct$ consisting of all fibrant objects homotopy orthogonal
to $\cm$. Since $\map(B,A)$ is weakly equivalent to $\hom(B,A)$ whenever $B$ is cofibrant and $A$ is fibrant 
and all morphisms from $\cm$ have cofibrant domains and codomains, we have
$$
\HAlg(\ct)=\Pre(\ct^{\op})\cap\cm^\bot.
$$
By \cite{CC}, 1.1, there is a functor $L:\SSet^\ct\to\cm^\bot$ preserving weak equivalences and 
equipped with a simplicial natural transformation $\eta:\Id\to L$ which is idempotent up to homotopy and, moreover, it is
a pointwise $\cm$-equivalence. The latter means that each object from $\cm^\bot$ is homotopy orthogonal to $\eta_K$
for each $K$ in $\SSet^\ct$. 

Consider a diagram $D:\cd\to\HAlg(\ct)$ and a $\ct$-algebra $A$. We have (where $R_c$ denotes a cofibrant replacement
functor in $\SSet^\ct$)
\begin{align*}
\hom(R_cL(\hocolim_f D),A)&\simeq\map(R_cL(\hocolim_f D),A)\\
&\simeq\map(L(\hocolim_f D),A)\\
&\simeq\map(\hocolim_f D,A)\\
&\simeq\hom(\hocolim_f D,A)\\
&\simeq\holim_s\hom(D,A).
\end{align*}
Thus $R_cL(\hocolim_f D)$ is a fibrant homotopy colimit of $D$ in $\HAlg(\ct)$. Hence $\HAlg(\ct)$ has fibrant
homotopy colimits. By \ref{th3.5}, we get a simplicial functor 
$$
F:\Pre(\ct^{\op})\to\HAlg(\ct)
$$
which preserves fibrant homotopy colimits and whose composition with $Y_{\ct^{\op}}$ is homotopy equivalent to the codomain
restriction of $Y_{\ct^{\op}}$ to $\HAlg(\ct)$. Consequently, $F$ is homotopy left adjoint to the inclusion
$G:\HAlg(\ct)\to\Pre(\ct^{\op})$. It follows from \ref{prop4.9} and \ref{prop4.14} that $F$ preserves homotopy strongly 
finitely presentable objects. Since, by \ref{prop4.11},  every object of $\Pre(\ct^{\op})$ is a homotopy sifted fibrant
homotopy colimit of homotopy strongly finitely presentable objects, $\HAlg(\ct)$ has the same property (because $F$ preserves 
fibrant homotopy colimits by \ref{lem3.3}). Hence $\HAlg(\ct)$ is a homotopy variety.

II. Let $\ck$ be a homotopy variety and $\ca$ be a set from \ref{def4.10}. Let $\bar\ca$ be the closure of $\ca$ under 
finite coproducts in $\ck$. By \ref{prop4.8}, each object of $\bar\ca$ is homotopy strongly finitely presentable in $\ck$. 
Put $\ct=(\bar\ca)^{op}$. Then $\ct$ is a simplicial algebraic theory. Let
$$
E:\ck\to\Pre(\bar\ca)
$$
be the simplicial functor given by
$$
E(K)=\hom(-,K)
$$
where the hom-functor is restricted to $\bar\ca$. Clearly, $E$ has values in $\HAlg(\ct)$ and we will denote its codomain 
restriction $\ck\to\HAlg(\ct)$ by $E$ as well. We will show that $E$ is an equivalence.

Consider objects $K_1$ and $K_2$ from $\ck$ and express them as homotopy sifted fibrant homotopy colimits
$K_i=\hocolim_f D_i$ of $D_i:\cd_i\to\ca$ where $i=1,2$. Then we have
\begin{align*}
\hom(K_1,K_2)&\simeq\hom(\hocolim_f D_1,\hocolim_f D_2)\\
&\simeq \holim_s\hom(D_1,\hocolim_f D_2)\\
&\simeq\holim_s\hocolim_f \hom(D_1,D_2)\\
&\simeq \holim_s\hocolim_f \hom(\hom(-,D_1),\hom(-,D_2))\\
&\simeq\holim_s \hom(\hom(-,D_1),\hocolim_f \hom(-,D_2))\\
&\simeq \hom(\hocolim_f \hom(-,D_1),\hocolim_f \hom(-,D_2))\\
&\simeq \hom(\hom(-,\hocolim_f D_1),\hom(-,\hocolim_f D_2))\\
&\simeq \hom(EK_1,EK_2).
\end{align*}
Here, we have used the homotopy invariance of simplicial homotopy colimits, the enriched Yoneda lemma, the homotopy absolute 
presentability of hom-functors in $\Pre(\bar\ca)$ (see \ref{re3.6}) and homotopy strong finite presentability of objects from 
$\ca$. Hence $E$ satisfies the first condition in the definition of an equivalence. The second condition follows from 
the fact that $\Pre(\bar\ca)$ is the free completion of $\bar\ca$ under fibrant homotopy colimits (see \ref{th3.5}) because 
it provides a simplicial functor $F:\Pre(\bar\ca)\to\ck$ with $FE\cong\Id_\ck$. Then its domain restriction $F:\HAlg(\ct)\to\ck$
satisfies $FE\cong\Id_\ck$, which yields the second condition.
\end{proof} 
 
\begin{defi}\label{def4.16}
{\em

Let $\cc$ be a small fibrant simplicial category. Then $\HSind(\cc)$ will denote the full subcategory of $\Pre(\cc)$ 
consisting of homotopy sifted fibrant homotopy colimits of hom-functors.  
}
\end{defi} 

\begin{theo}\label{th4.17}
Let $\cc$ be a small fibrant simplicial category having finite coproducts. Then the simplicial categories $\HAlg(\cc^{op})$ 
and $\HSind(\cc)$ are equivalent. 
\end{theo}
\begin{proof}
Since homotopy sifted fibrant homotopy colimits commute with finite products in $\Sp$, we always have
$$
\HSind(\cc)\subseteq\HAlg(\cc^{op}).
$$
Conversely, we know that each object from $\HAlg(\cc^{\op})$ is a homotopy sifted fibrant homotopy colimit of finite
coproducts of hom-functors (see the proof of \ref{th4.15}). Since $L(m_{X_1\dots X_n})$ is a weak equivalence for each
morphism $m_{X_1\dots X_n}$ from this proof (see \cite{DF}, 1.C.5), each object from $\HAlg(\cc^{\op})$ is homotopy 
equivalent to an object from $\HSind(\cc)$.
\end{proof}

\begin{rem}\label{re4.18}
{\em
As a consequence, we get that $\HSind(\cc)$ has all homotopy sifted fibrant homotopy colimits. Thus it
is the free completion of $\cc$ under homotopy sifted fibrant homotopy colimits for each small simplicial category
$\cc$ having finite coproducts. Hence $\HSind(\cc)$ is analogous to the free completion $\Sind(\cc)$ of a category $\cc$ 
under sifted colimits introduced in \cite{AR1}.
}
\end{rem}

\section{Homotopy locally finitely presentable categories}

\begin{defi}\label{def5.1}
{\em
A fibrant simplicial category $\ck$ will be called \textit{homotopy locally finitely presentable} provided that it has 
fibrant homotopy colimits and has a set $\ca$ of homotopy finitely presentable objects such that every object of $\ck$ 
is a filtered fibrant homotopy colimit of objects from $\ca$. 
}
\end{defi} 

\begin{propo}\label{prop5.2}
Every homotopy variety is homotopy locally finitely presentable.
\end{propo}
\begin{proof}
Let $\ca$ be the set of homotopy strongly finitely presentable objects from \ref{def5.1} and $\bar\ca$ be
its closure under finite fibrant homotopy colimits in $\ck$. By \ref{re4.7} and \ref{prop3.9}, $\bar\ca$
consists of homotopy finitely presentable objects. Since each object in $\ck$ is a fibrant homotopy colimit of 
objects from $\ca$ and a fibrant homotopy colimit can be expressed as a filtered fibrant homotopy colimit of finite 
fibrant homotopy colimits, each object of $\ck$ is a filtered fibrant homotopy colimit of objects from $\bar\ca$. Hence 
$\ck$ is homotopy locally finitely presentable.
\end{proof}

\begin{coro}\label{cor5.3}
Let $\cc$ be a small fibrant simplicial category. Then the category $\Pre(\cc)$ is homotopy locally finitely presentable. 
\end{coro}
\begin{proof}
It follows from \ref{prop4.11} and \ref{prop5.2}. 
\end{proof}

\begin{defi}\label{def5.4}
{\em
A \textit{finite homotopy limit theory} is defined as a small fibrant simplicial category $\ct$ having all finite fibrant
homotopy limits. 

A \textit{homotopy} $\ct$-\textit{model} is a simplicial functor $A:\ct\to\Sp$ belonging to $\Pre(\ct^{\op})$ and
preserving finite fibrant homotopy limits.

We will denote by $\HMod(\ct)$ the full subcategory of $\Pre(\ct^{\op})$ consisting of all homotopy $\ct$-models.
}
\end{defi} 

\begin{propo}\label{prop5.5}
Let $\ct$ be a finite homotopy limit theory. Then the simplicial category $\HMod(\ct)$ is closed in $\Pre(\ct^{\op})$ both 
under fibrant homotopy limits and filtered fibrant homotopy colimits.
\end{propo}
\begin{proof}
It is analogous to that of \ref{prop4.14} (using \ref{prop3.8}). 
\end{proof}

\begin{theo}\label{th5.6}
A fibrant simplicial category $\ck$ is homotopy locally finitely presentable iff it is equivalent to $\HMod(\ct)$ 
for some finite homotopy limit theory $\ct$.
\end{theo}
\begin{proof}
I. Let $\ct$ be a finite homotopy limit theory. We proceed analogously as in the proof of \ref{th4.15}. We only need 
to replace the morphisms $m_{X_1\dots X_n}$ from that proof by morphisms
$$
m_D:\hocolim_f\hom(D,-)\to\hom(\holim_f D,-)
$$
for each finite diagram $\cd\to\ct$. By the dual of \ref{re3.1}(b), $m_D$ corresponds to the morphism
$$
\bar{m}_D:B((-\downarrow\cd^{\op})^{\op})\to\hom(\hom(D,-),\hom(\holim_fD,-)).
$$
Since the domain of $\bar{m}_D$ is isomorphic to $B(\cd\downarrow -)$ and the codomain to $\hom(\holim_f D,D)$,
$\bar{m}_D$ corresponds to the morphism
$$
\tilde{m}_D:B(\cd\downarrow -)\to\hom(\holim_f D,D).
$$
Now, in order to define $m_D$, we take the morphism $\tilde{\delta}_D$ from \ref{re2.4} for $\tilde{m}_D$.

II. Let $\ck$ be homotopy locally finitely presentable simplicial category and $\ca$ be the set from \ref{def5.1}.
Let $\bar\ca$ be the closure of $\ca$ under finite fibrant homotopy colimits in $\ck$. By \ref{prop3.9}, each
object from $\bar\ca$ is homotopy finitely presentable in $\ck$. Now, we put $\ct=(\bar\ca)^{\op}$ and proceed
analogously as in the proof of \ref{th4.15}.
\end{proof}

\begin{coro}\label{cor5.7}
A homotopy locally finitely presentable category has all fibrant homotopy limits.
\end{coro}
\begin{proof}
It follows from \ref{th5.6} and \ref{prop5.5}.
\end{proof}

\begin{rem}\label{re5.8}
{\em
By a \textit{homotopy finite limit sketch} is meant a triple $\ch=(\ct,{\bf L},\sigma)$ consisting of a small
fibrant simplicial category $\ct$, a set $\bf L$ of finite diagrams in $\ct$ and an assignment $\sigma$ of
a morphims
$$
\sigma(D):B(\cd\downarrow -)\to\hom(X_D,D)
$$
in $\SSet^{\cd}$ to each diagram $D\in\bf L$.

By a \textit{homotopy model} of $\ch$ is meant a simplicial functor $A:\ct\to\Sp$ belonging to $\Pre(\ct^{\op})$ and
sending $\sigma(D)$ to $\tilde{\delta}_D$ for each $D\in{\bf L}$.  

We will denote by $\HMod(\ch)$ the full subcategory of $\Pre(\ct^{\op})$ consisting of all homotopy models of $\ch$.

Every homotopy finite limit theory is a homotopy finite limit sketch. Since the part I. of the proof of \ref{th5.6} is
valid for each homotopy finite limit sketch $\ch$, $\HMod(\ch)$ is always homotopy locally finitely presentable.
}
\end{rem}

\begin{defi}\label{def5.9}
{\em
Let $\cc$ be a small fibrant simplicial category. Then $\HInd(\cc)$ will denote the full subcategory of 
$\Pre(\cc)$ consisting of filtered fibrant homotopy colimits of hom-functors. 
}
\end{defi} 

\begin{theo}\label{th5.10}
Let $\cc$ be a small fibrant simplicial category having finite fibrant homotopy colimits. Then the simplicial
categories $\HInd(\cc)$ and $\HMod(\cc^{op})$ are equivalent. 
\end{theo}
\begin{proof}
Since filtered fibrant homotopy colimits commute with finite fibrant homotopy limits in $\Sp$ (see \ref{prop3.8}), 
we always have
$$
\HInd(\cc)\subseteq\HMod(\cc^{op}).
$$
Conversely, we know that each object from $\HMod(\cc^{\op})$ is a filtered fibrant homotopy colimit of finite fibrant
homotopy colimits of hom-functors (using \ref{th3.5}). Since $L(m_D)$ is a weak equivalence for each
morphisms $m_D$ from the proof of \ref{th5.6} (see \cite{DF}, 1.C.5), each object from $\HMod(\cc^{\op})$ is homotopy 
equivalent to an object from $\HInd(\cc)$.
\end{proof}

\begin{rem}\label{re5.11}
{\em
(1) As a consequence, we get that $\HInd(\cc)$ has all filtered fibrant homotopy colimits. Thus it
is the free completion of $\cc$ under filtered fibrant homotopy colimits. Hence $\HInd(\cc)$ is
analogous to the free completion $\Ind(\cc)$ of a category $\cc$ under filtered colimits introduced in
\cite{AGV}.

(2) Everything in this section can be done for an arbitrary regular cardinal $\lambda$ instead of $\omega$.
It means that we work with homotopy $\lambda$-filtered fibrant homotopy colimits and compare homotopy
locally $\lambda$-presentable categories with categories of models of $\lambda$-small homotopy limit
theories.
}
\end{rem}

\end{document}